\documentclass[12pt]{amsart}

\def\colon{{:}\;}
\newtheorem{theo}{Theorem}
\newtheorem{prop}[theo]{Proposition}

\newtheorem{coro}[theo]{Corollary}

\theoremstyle{definition}

\newtheorem{exam}[theo]{Example}

\theoremstyle{remark}

\newcommand{\ZZ}{\mathbb Z}

\newcommand{\MM}{\mathbb M}

\newcommand{\bT}{\mathbb T}

\newcommand{\zd}{\mathbb Z^d}

\newcommand{\cA}{\mathcal A}

\newcommand{\cP}{\mathcal P}

\newcommand{\n}{{\bf n}}
\newcommand{\m}{{\bf m}}

\newcommand{\e}{{\bf e}}

\newcommand{\bo}{{\bf 0}}

\newcommand{\h}{\operatorname{h}}

\def\eref(#1){\textup{(\ref{#1})}}

\begin{document}\allowdisplaybreaks\frenchspacing

\title[Isomorphism rigidity
in entropy rank two]{Isomorphism rigidity\\in entropy rank two}
\author{Manfred Einsiedler}
\author{Thomas Ward}
\date{\today}
\thanks{The second named author thanks the University of
Washington for their hospitality and the support of NSF grant DMS
0222452} 

\begin{abstract}
We study the rigidity properties of a class of algebraic $\mathbb
Z^3$-actions with entropy rank two. For this class, conditions are
found which force an invariant measure to be the Haar measure on
an affine subset. This is applied to show isomorphism rigidity for
such actions, and to provide examples of non-isomorphic
$\ZZ^3$-actions with all their $\mathbb Z^2$-sub-actions
isomorphic. The proofs use lexicographic half-space entropies
and total ergodicity along critical directions.
\end{abstract}

\maketitle

\section{Introduction}

An {\bf algebraic $\zd$-action} is a $\zd$-action on a compact
abelian metrizable group by automorphisms. Rigidity for such
actions is a circle of results that give explicit descriptions of
all invariant measures for $\alpha$, or all measurable
isomorphisms between such systems, under certain hypotheses for $d>1$. The
case of a $\zd$-action by toral automorphism is well studied
(see~\cite{MR2003h:37007} for isomorphism rigidity
and~\cite{Einsiedler-Lindenstrauss-RA},
\cite{MR2002i:37035}, \cite{MR97d:58116} for measure rigidity); in the
toral case individual elements of the action have finite entropy.
For actions on zero-dimensional groups the case of a general
irreducible action was first studied in~\cite{MR2001j:37004},
where again the individual elements of the action have finite
entropy. The general case of mixing algebraic $\zd$-actions on
zero-dimensional groups with zero entropy was studied
in~\cite{Bhattacharya02}, \cite{BhattacharyaSchmidt02}
and~\cite{Einsiedler02}. Another type of rigidity --- differences
between apparently related zero-dimensional systems forcing them
to be disjoint
--- is studied in~\cite{EinsiedlerWarddisjoint} also using entropy
methods.

Our purpose is to study isomorphism rigidity for a particular
class of algebraic $\ZZ^3$-actions on connected groups with the
property that all the $\ZZ^2$-sub-actions have positive entropy.
We follow the path taken in~\cite{Einsiedler02}, and avoid repetition
of certain technicalities by citing results from that paper.
The
special class of systems studied allows a transparent proof,
but it is clear that the the underlying rigidity phenomena is
more extensive than what is shown here suggests.

The systems we study, $ET$ systems, have two
distinguished polynomial parameters.
The first of these, $g$, is
$E$xpanding, and the second, $f$, is $T$riangular.
The assumption
on $g$ is essential for the method used here, while the assumption
on $f$ does not seem to be essential but does significantly simplify
the arguments.
The paradigmatic example of an $ET$ system is the
so-called space helmet~(cf.~\cite[Example~5.8]{MR1869066};
Example~\ref{Nexample}).

Let $\mathsf X_i=(X_i,\alpha_i)$ be an algebraic $\ZZ^d$-action
for $i=1,2$.
A factor map $\varphi\colon\mathsf X_1\rightarrow\mathsf
X_2$ is a (Borel) measurable map from $X_1$ onto $X_2$ with
$\varphi\circ\alpha_1^\n(x)=\alpha_2^\n\circ\varphi(x)$ for a.e.\
$x\in X_1$ and all $\n\in\ZZ^d$.
A factor map $\varphi\colon\mathsf
X_1\rightarrow\mathsf X_2$ is\begin{itemize}
    \item a {\bf conjugacy} if it is invertible;
    \item {\bf affine} if $\varphi(x)=\varphi_g(x)+y$ for
    some $y\in X_2$ and continuous group homomorphism
$\varphi_g\colon X_1\rightarrow X_2$; and is
    \item an {\bf algebraic isomorphism} if it is an isomorphism
    of the groups $X_1$ and $X_2$.
\end{itemize}
We will always use $\lambda_X$ to denote the Haar measure on a compact
abelian group $X$.

A class of systems exhibits {\bf isomorphism rigidity} if every
conjugacy or factor map is affine, and exhibits {\bf measure
rigidity} if very mild additional assumptions on an
invariant ergodic Borel
measure (positive entropy, for example) force it to be a translate of the
Haar measure of a closed
subgroup. 

We use a standard approach to the description of algebraic
$\ZZ^d$-actions (see~\cite{MR97c:28041} for more background).
Let
$R_d=\ZZ[u^{\pm 1}_1,\dots,u^{\pm 1}_d]$ be the ring of Laurent
polynomials in $d$ commuting variables with integer coefficients,
and write ${\mathbf u}^{\n}$ for the monomial $u_1^{n_1}\dots
u_d^{n_d}$.
A polynomial $f\in R_d$ is a sum
\[
  f(\mathbf u)=\sum_{\n\in\zd}f_\n{\mathbf u}^\n,
\]
where $f_\n\in\ZZ$ for all $\n\in\zd$ and $f_\n=0$ for all but
finitely many $\n$.
Any $R_d$-module $M$ defines an algebraic
$\mathbb Z^d$-action $\mathsf X_M=(X_M,\alpha_M)$ as follows.
The compact group
$X_M=\widehat{M}$ is the (Pontryagin) character group of the
additive group $M$; the action $\alpha_M$ is defined by
$\alpha_M^{\n}(x)=\widehat{\beta^{\n}}(x)$ where
$\beta^{\n}(m)={\mathbf u}^{\n}m$.
By duality, any algebraic
$\ZZ^d$-action $\mathsf X$ similarly defines an $R_d$-module
$M_{\mathsf X}$; the module structure is given by defining the
product of $f\in R_d$ with $a\in M_{\mathsf X}$ to be
\[
 fa=\sum_{\n\in\zd}f_\n\widehat{\alpha^\n} a.
\]
This gives a one-to-one correspondence between algebraic
$\ZZ^d$-actions on compact (metrizable) abelian groups and
(countable) $R_d$-modules. 

A system $\mathsf X_M$ is called {\bf prime} if $M$ is a cyclic
$R_d$-module of the form $R_d/P$ for some prime ideal $P\subset
R_d$.
Notice that by duality, prime systems $\mathsf X_M$ and
$\mathsf X_N$ are algebraically isomorphic if and only if their
defining ideals are equal.

\begin{exam}\label{Mexample} The module $M=R_2/\langle 1+u_1+u_2\rangle$
defines an algebraic $\ZZ^2$-action as follows.
The dual of $M$ is
\[
X_M=\bigl\{x\in\bT^{\ZZ^2}\colon
x_\n+x_{\n+(1,0)}+x_{\n+(0,1)}=0\mbox{
for any }\n\in\ZZ^2\bigr\},
\]
and the $\ZZ^2$-action $\alpha_M$ is the restriction of the shift
action
\[
 \alpha_M^\n(x)^{\vphantom \n}_\m=
x^{\vphantom \n}_{\m+\n} \mbox{ for any }\n,\m\in\ZZ^2
\]
to $X_M$.
By~\cite{MR96d:22004} (see
also~\cite[Th.~20.8, Th.~23.1]{MR97c:28041},
\cite{MR93h:28032})
$\alpha$ is isomorphic to a Bernoulli shift,
so there are many non-affine conjugacies from $\mathsf
X_M$ to itself: the system does not exhibit isomorphism rigidity.
\end{exam}

\begin{exam}\label{Nexample} The module $N=R_3/\langle
1+u_1+u_2,u_3-2\rangle$ defines an algebraic action as follows.
Let
\[
X_N=\bigl\{x\in\bT^{\ZZ^3}\colon
x_\n+x_{\n+\e_1}+x_{\n+\e_2}=0,x_{\n+\e_3}=2x_\n
\mbox{ for }\n\in\ZZ^3\bigr\},
\]
and define the $\ZZ^3$-action $\alpha_N$ to be the restriction of
the usual shift action to $X_N$.
This example and the behaviour of
its lower-rank sub-actions was studied
in~\cite[Ex.~5.8]{MR1869066}; the structure of its non-expansive
subdynamics has earned $\mathsf X_N$ the sobriquet of `space
helmet'.
\end{exam}
The $\ZZ^3$-action $\alpha_N$ in Example~\ref{Nexample} can be
obtained from the $\ZZ^2$-action
$\alpha_M$ in Example~\ref{Mexample} by defining the
third generator to be multiplication by $2$ on each
coordinate, and then passing to
the invertible extension of this non-invertible map.
Thus a
conjugacy from $\mathsf X_M$ to itself extends to a conjugacy from
$\mathsf X_N$ to itself if and only if it commutes with
multiplication by $2$.
Theorem~\ref{t:main} shows that
very few maps can have this property;
in particular it implies the following proposition.

\begin{prop}\label{p:space}
Every conjugacy from $\mathsf X_N$ to itself is affine.
\end{prop}

In the next section we define the class $ET$ and show how
isomorphism rigidity may be deduced from the technical
Proposition~\ref{p:meas}.
Section~\ref{s:examples} gives examples
of algebraically non-isomorphic $\ZZ^3$-actions for which all
$\ZZ^2$-sub-actions are conjugate.
Section~\ref{s:proof} contains
the proof of Proposition~\ref{p:meas}, following the ideas
in~\cite{Einsiedler02}.

\section{The class $ET$ and isomorphism rigidity}\label{s:class}

An algebraic $\ZZ^3$-action $\mathsf X_M$ is said to be of class
$ET$ if $M$ is Noetherian and there are
polynomials $f\in\mathbb Z[u_1,u_2]$
and $g\in\mathbb Z[u_3]$
with the following properties.
\begin{itemize}
\item Both $f$ and $g$ annihilate $M$. \item $g$ is monic, and
every zero $z$ of $g$ has $\vert z\vert>1$ ($g$ is $E$xpanding).
\item The Newton polygon of $f$ (the convex hull of $\{\mathbf
n\colon f_{\mathbf n}\neq0\}$) is a triangle with corners at
$(0,0,0)$, $(a,0,0)$, $(0,a,0)$ for some $a>0$ ($f$ is
$T$riangular) and the coefficients corresponding to these corners
are $\pm 1$. 
\end{itemize}
Notice that property $ET$ is stable in the following
sense.
If $\mathsf X$ and $\mathsf Y$ are $ET$, then so is $\mathsf X\times \mathsf
Y$.
If $M$ is a Noetherian module with the property that $\mathsf
X_{R_d/P}$ is $ET$ for every prime ideal associated to $M$, then
$\mathsf X_M$ is $ET$.
In both cases the reason is that the
product of two expanding polynomials is expanding, and the
product of two triangular polynomials is triangular.

\begin{theo}\label{t:main}
 Let $\mathsf X_1$ and $\mathsf X_2$ be $ET$ algebraic $\ZZ^3$-actions.
 Suppose $\alpha_1$ is mixing and the $\ZZ^2$-sub-action generated by
 $\alpha_1^{\e_2}$ and $\alpha_1^{\e_3}$ has completely positive
 entropy.
Then every factor map
 from $\mathsf X_1$ to $\mathsf X_2$ is affine.
\end{theo}

\begin{coro}\label{c:prime}
 Let $P_1\neq P_2$ be prime ideals in $R_3$.
Assume that
$\mathsf X_{R_3/P_2}$ is $ET$,
$\mathsf
 X_{R_3/P_1}$ is $ET$ and mixing, and the $\ZZ^2$-sub-action generated by
 $\alpha_{R_3/P_1}^{\e_2}$ and $\alpha_{R_3/P_1}^{\e_3}$
 has completely positive entropy.
Then $\mathsf X_{R_3/P_1}$ and $\mathsf X_{R_3/P_2}$ are not
 measurably conjugate.
\end{coro}

\begin{proof}  By Theorem~\ref{t:main}, it is enough to show that
 the two actions are not algebraically isomorphic.
 If $\varphi$ were an algebraic
 isomorphism between $\mathsf X_{R_3/P_1}$ and $\mathsf X_{R_3/P_2}$,
 then the dual map $\widehat{\varphi}\colon
R_3/P_2\rightarrow R_3/P_1$
 would be a module isomorphism, which would imply that $P_1=P_2$.
\end{proof}

Theorem~\ref{t:main}
is proved using the restricted version of measure
rigidity stated in Proposition~\ref{p:meas}.
An
$\alpha_M$-invariant measure $\mu$ is {\bf totally ergodic} if
$\alpha_M^\n$ is an ergodic transformation of $(X_M,\mu)$ for all
$\n$ in $\ZZ^3\backslash\{0\}$.

\begin{prop}\label{p:meas}
 Let $\mathsf X_M$ be an $ET$ algebraic $\ZZ^3$-action.
 Suppose $\mu$ is a totally ergodic $\alpha_M$-invariant measure.
 Then there exists a closed $\alpha_M$-invariant subgroup $Y\subset X_M$
 such that $\mu$ is invariant under translation by elements of
 $Y$, and the action generated by $\alpha_M^{\e_2}$ and
 $\alpha_M^{\e_3}$ induced on the factor $X_M/Y$ has zero entropy with respect
 to $\mu$.
\end{prop}

Theorem~\ref{t:main} follows by a well-known argument
of Thouvenot from
Proposition~\ref{p:meas} (see~\cite{Einsiedler02}
or~\cite{MR2003h:37007} for more details).
A sketch of this
argument follows.
As noted above,
if $\mathsf X_1$ and $\mathsf X_2$ are $ET$, then so is
$\mathsf X_1\times\mathsf X_2$.
A conjugacy $\varphi$ from
$\mathsf X_1$ to $\mathsf X_2$ defines a joining $\mu$ supported
on the graph $G_{\varphi}$ of $\varphi$.
The projection map
$\pi_1\colon X_1\times X_2\rightarrow X_1$ onto the first coordinate
satisfies $\mu(A)=\lambda_{X_1}(\pi_1(A\cap G_\phi))$ for any
$A\subset X_1\times X_2$.
The measure~$\mu$ satisfies the
assumption of Proposition~\ref{p:meas} since $\alpha_1$ is mixing.
By assumption, the sub-action generated by $\alpha_1^{\e_2}$ and
$\alpha_1^{\e_3}$ has completely positive entropy.
By
Proposition~\ref{p:meas}, the entropy on the factor $X_1\times
X_2/Y$ vanishes, so this factor must be trivial.
It follows
that~$\mu$ is the Haar measure of the affine subset $x_0+Y$ for
some $x_0\in X_1\times X_2$, and the graph $G_\varphi$ agrees with
the affine subset $x_0+Y$ a.e.
We conclude that $\varphi$ is
affine $\lambda_{X_1}$-a.e.

\section{Examples}\label{s:examples}

We give examples of algebraically --- hence, by
Theorem~\ref{t:main}, measurably --- non-isomorphic
$\ZZ^3$-actions all of whose $\ZZ^2$-sub-actions are measurably
conjugate.

To study the action of a subgroup $\Lambda\subset\ZZ^3$, it is
useful to define the corresponding sub-ring
 \[
  R_\Lambda=\ZZ[u^\n\colon\n\in\Lambda]\subset R_3.
 \]
Dynamical properties of the sub-action $\alpha_M\vert_\Lambda$ are
governed by the structure of the module $M$ over the ring
$R_\Lambda$.
Notice that in the case of a prime system $\mathsf
X_{R_d/P}$, the only prime ideal in $R_{\Lambda}$ associated to
$R_d/P$ as a module over $R_\Lambda$ is $P\cap R_\Lambda$.

\begin{exam}\label{e:plusminus2}
 Let
$$
P_1=\langle 1+u_1+u_2, u_3-2\rangle
\mbox{ and }
P_2=\langle 1+u_1+u_2,
 u_3+2\rangle.$$
We will show that the prime actions
 $\mathsf X_{R_3/P_i}$ for $i=1,2$ are mixing and have zero entropy.
 Furthermore, the sub-actions corresponding to subgroups
 of rank two have completely positive entropy,
 are isomorphic to Bernoulli shifts and are conjugate.
 By Corollary~\ref{c:prime}, the two actions are not conjugate.

 The ideals are prime ideals
 since in each case
$$
R_3/P_i=\ZZ[1/2,u_1^{\pm 1},u_2^{\pm 1}]/\langle 1+u_1+u_2\rangle.
$$
 The only difference between the two actions is that
 $\alpha_{R_3/P_1}^{\e_3}$ acts by multiplication by $2$ and
 $\alpha_{R_3/P_2}^{\e_3}$ acts by multiplication by $-2$.
This
 shows that $\alpha_{R_3/P_1}^{2\e_3}$ and
 $\alpha_{R_3/P_1}^{2\e_3}$ act identically, so
 the restrictions of the actions to the subgroup
 $\ZZ^2\times(2\ZZ)$ are algebraically isomorphic, and the entropies
 for any rank two subgroup of $\ZZ^2\times(2\ZZ)$ must therefore
 agree.
 Since the entropy of the $\Lambda$-sub-action and the entropy
 for the $2\Lambda$-sub-action determine each other,
 it follows that the
 sub-actions $\alpha_{R_3/P_1}\vert_\Lambda,\alpha_{R_3/P_1}\vert_\Lambda$
 induced by any subgroup $\Lambda\subset\ZZ^3$ of rank two have equal entropy.

 We claim that both $\mathsf X_{R_3/P_1}$ and $\mathsf X_{R_3/P_2}$
are mixing and have zero entropy.
 For mixing, we need to show that $u^\n-1\notin P_i$ for every
 $\n$ in $\ZZ^3\setminus\{\bo\}$.
Assume that $u^\n-1\in P_i$, and let
 $\zeta_3=e^{2\pi i/3}$.
Then $(\zeta_3,\zeta_3^2, \pm 2)$
 belongs to the variety of $P_i$.
This shows that
 $\zeta_3^{n_1+2n_2}(\pm 2)^{n_3}=1$, so $n_3=0$.
Going through the
 argument again using the points $(1,-2,\pm 2),(-2,1,\pm 2)$ on the variety shows that
 $\n=0$.
The action has zero entropy since it is a prime action
 whose prime ideal is not principal (see~\cite{MR92j:22013}
 or~\cite[Cor.~18.5]{MR97c:28041}).

 Let $\Lambda\subset\ZZ^3$ be any subgroup isomorphic
 to $\ZZ^2$.
We claim that the sub-action
 $\alpha_{R_3/P_i}\vert_\Lambda$ is isomorphic to a two-dimensional
 Bernoulli shift for $i=1,2$.
 By~\cite{MR96d:22004} (see also~\cite[Th.~23.1]{MR97c:28041}),
 an algebraic $\zd$-action is isomorphic to a Bernoulli shift if and only if it has
 completely positive entropy.
By~\cite{MR92j:22013},
 (see also~\cite[Th. 20.8]{MR97c:28041}) this is the case if and
 only if every associated prime ideal is generated by a single
 polynomial which is not cyclotomic.
Applying this to the
 sub-action $\alpha_{R_3/P_i}\vert_\Lambda$, it is enough to show that
 $P_i\cap R_\Lambda$ is generated by a single polynomial which is
 not cyclotomic.

 Assume first that $P_i\cap R_\Lambda$ is trivial.
In this case
 $R_\Lambda$ is a sub-ring of $R_3/P_i$.
This is impossible
 since $R_\Lambda$ has
 transcendence degree two but $R_3/P_i$ has transcendence degree one,
 so $P_i\cap R_\Lambda$ cannot be trivial.
Similar arguments
 show that $P_i\cap R_\Lambda$ must be a principal ideal.
By the
 above we already know that $P_i$ does not contain any
 cyclotomic polynomial, so $P_i\cap R_\Lambda$ is not generated
 by a cyclotomic polynomial.
\end{exam}

In the remaining examples we will not repeat the arguments above
showing that the action is a prime action, has zero entropy, is
mixing and that every sub-action for a rank two subgroup is
Bernoulli.
In each case these use well-known algebraic methods and
characterizations from~\cite{MR97c:28041}.
We will just check the
following property: systems $\mathsf X_1$ and $\mathsf X_2$ are
{\bf $\ZZ^2$-entropy equivalent} if for every subgroup
$\Lambda\subset\ZZ^3$ of rank two,
$\h(\alpha_1\vert_\Lambda)=\h(\alpha_2\vert_\Lambda)$.

\begin{exam}\label{e:rescaled}
 Let
$$
M=R_3/\langle 1+u_1+u_2, u_3-4\rangle
\mbox{ and }
N=R_3/\langle 1+u_1^2+u_2^2,u_3-2\rangle.
$$
As before, $\mathsf X_M$ and $\mathsf X_N$
 are mixing, zero entropy prime actions with the property that every sub-action
 of a rank two subgroup has completely positive entropy.
 By Corollary~\ref{c:prime}, the two actions are not conjugate.

 We claim the actions are $\ZZ^2$-entropy equivalent.
To see this,
 consider the sub-action $\alpha_N\vert_\Gamma$ for
 $\Gamma=(2\ZZ)^3$ and use the rescaled variables $w_j=u_j^2$ and
 the ring $S_3=\ZZ[w_1^{\pm 1},w_2^{\pm 1},w_3^{\pm 1}]$ to study the sub-action
 (see~\cite{MR98d:58053} for other applications of rescaling).
 Define a homomorphism of $S_3$-modules by
 \[
  \varphi\colon \bigl(S_3/\langle
  1+w_1+w_2,w_3-4\rangle\bigr)^4\longrightarrow
  N
 \]
 by
 \[
 \varphi(a_1,a_2,a_3,a_4)=a_1+u_1a_2+u_2a_3+u_1u_2a_4.
 \]
 This is an isomorphism, which shows
 that the action $\alpha_N\vert_\Gamma$ is algebraically
 isomorphic to four disjoint copies of $\alpha_M$.
Now
let $\Lambda\subset\ZZ^3$ be a subgroup isomorphic to $\ZZ^2$.
 Clearly $2\Lambda\subset\Gamma$, and the algebraic isomorphism dual to $\varphi$ carries
 $\alpha_N\vert_{2\Lambda}$ to $\alpha_M\vert_\Lambda$.
 It follows that
 \[
  \h (\alpha_N\vert_\Lambda)=\textstyle\frac{1}{4}\h (\alpha_N\vert_{2\Lambda})
   =\h (\alpha_M\vert_\Lambda),
 \]
 showing that $\mathsf X_M$ and $\mathsf X_N$ are $\ZZ^2$-entropy
 equivalent.
\end{exam}

In Examples~\ref{e:plusminus2} and~\ref{e:rescaled} the pairs of
ideals were closely related.
In the following the same triangular
polynomial $f(u_1,u_2)=1+u_1+u_2$ is used, but the expanding
polynomials are chosen with rationally independent roots.

By~\cite{MR92j:22013} or~\cite[Th.~18.1]{MR97c:28041}, for any
non-zero polynomial $h\in R_d$ the entropy of the action
$\alpha_{R_d/\langle h\rangle}$ is given by
\begin{equation*}
 \h(\alpha_{R_d/\langle h\rangle})=\log\MM(h),
\end{equation*}
where the Mahler measure $\MM(h)$ of the polynomial $h$ is defined
by
\begin{equation*}
 \log\MM(h)=\int_0^1\!\!\cdots\!\int_0^1\log|h(e^{2\pi it_1},\ldots,e^{2\pi
 it_d})|\ \operatorname{d}\!t_1\cdots \operatorname{d}\!t_d.
\end{equation*}
For $d=1$ the Mahler measure can be directly expressed in terms of
expanding eigenvalues.
Jensen's formula (see~\cite[Lemma~1.8]{MR2000e:11087}
or~\cite[Prop.~16.1]{MR97c:28041})
shows that
\begin{equation}\label{e: Jensen}
 \int_0^1\log |h(e^{2\pi i
 t})|\operatorname{d}\!t=\sum_{j=1}^s\log^+|\zeta_j|+\log|a_s|
\end{equation}
for a polynomial $h(u)=a_s\prod_{i=1}^{s}(u-\zeta_i)$.
This shows
that for $d=1$ the Mahler measure of $h$ only depends on the
absolute values of the zeros of $h$ and its leading coefficient.

\begin{exam}
 Let $g_1(u_3)=u_3^2+2u_3+10$ and $g_2(u_3)=u_3^2+4u_3+10$.
Eisenstein's
 criterion for $p=2$ shows that both polynomials are irreducible
 in $\ZZ[u_3]$.
If
 $g_i(u_3)=(u_3-\zeta_{i,1})(u_3-\zeta_{i,2})$ then
 $\zeta_{i,j}\notin\mathbb R$ and
 $\vert\zeta_{i,j}\vert^2=10$ for $i,j=1,2$.

 Let $M_1=R_3/\langle 1+u_1+u_2,g_1(u_3)\rangle$
and $M_2=R_3/\langle 1+u_1+u_2,g_2(u_3)\rangle$.
 As before $\mathsf X_{M_i}$ is a prime, mixing, zero
 entropy $ET$ system for which every sub-action
 of a rank two subgroup has completely positive entropy for $i=1,2$.
 By Corollary~\ref{c:prime}, the two actions are not conjugate.

 We show that $\mathsf X_{M_1}$ and $\mathsf X_{M_2}$ are
 $\ZZ^2$-entropy equivalent, so that the sub-actions for any subgroup of rank two
 are conjugate.
For $\Lambda=\ \ll\!\!\e_1,\e_2\!\!\gg$, this is clear
 since the associated prime ideal for this particular sub-action
 of $\mathsf X_{M_i}$ is
$$
\langle 1+u_1+u_2,g_i(u_3)\rangle\cap R_2=\langle 1+u_1+u_2\rangle.
$$
For
 a subgroup of rank two with $\mathbf e_3$ as a generator,
 the associated prime ideals for the sub-action are generated by
 $g_1$ and $g_2$, which are monic polynomials with two roots
 of absolute value $\sqrt{10}$.
Therefore $\MM(g_1)=\MM(g_2)=10$,
 and the entropy is a multiple of $\log 10$.
The multiplicative factor depends
 only on the geometry of $\Lambda$ and not on the polynomial $g_i$
 (see the general case below for more details).
 For a subgroup in general position, we need to argue that the
 Mahler measure of two
(in general unknown) polynomials coincide, and that
 two (potentially inscrutable) multiplicities coincide.

 To see what is creating the exact value of the entropy for a rank
 two subgroup, first consider the subgroup $\Lambda=\
 \ll\!\!\e_1,\e_2-\e_3\!\!\gg$.
Write $v_1=u_1$ and
 $v_2=u_2u_3^{-1}$.
The structure of the subsystem is given by the
 structure of $\widehat{\mathsf X_{M_i}}$ as a
 module over the ring $S=\ZZ[v_1^{\pm1},v_2^{\pm1}]$.
Let $c_1=-2$,
 $c_2=-4$ be the different traces in $g_1$, $g_2$.
 Since
$$
1+u_1+u_2=1+v_1+v_2u_3,
$$
the relation
\begin{equation}\label{e:tiltedfirstex}
(1\negthinspace +\negthinspace v_1^{\vphantom2}
\negthinspace +\negthinspace v_2^{\vphantom2}\zeta_{i,1}^{\vphantom2})
(1\negthinspace +\negthinspace v_1^{\vphantom2}
\negthinspace +\negthinspace v_2^{\vphantom2}\zeta_{i,2}^{\vphantom2})=
(1\negthinspace+\negthinspace v_1^{\vphantom2})^2+10v_2^2+
c_i^{\vphantom2}v_1^{\vphantom2}v_2^{\vphantom2}+
c_i^{\vphantom2}v_2^{\vphantom2}=0
\end{equation}
 holds for $i=1,2$.
Since this relation is irreducible over $\ZZ$, the
polynomial
in~\eqref{e:tiltedfirstex} generates
the only principal prime ideal associated to
 $M_i$ over $S$.
This means the whole subsystem is an
 invertible extension of the system defined by the
 relation~\eqref{e:tiltedfirstex} in the plane spanned by $\Lambda$.
 Thus
 \begin{eqnarray*}
  \h(\alpha_{M_i}\vert_{\Lambda})&=&\log\MM((1+v_1^{\vphantom2})^2
  +10v_2^2+c_i^{\vphantom2}v_1^{\vphantom2}v_2^{\vphantom2}
  +c_i^{\vphantom2}v_2^{\vphantom2}) \\
  &=&\log\MM(1+v_1+v_2\zeta_{i,1})+\log\MM(1+v_1+v_2\zeta_{i,2})\\
  &=&\log\vert\zeta_{i,1}\vert+\int_0^1\log^{+}\vert(1+e^{2\pi it})/\zeta_{i,1}\vert\mbox{d}t\\
  &&\quad\quad\quad\quad\quad+
  \log\vert\zeta_{i,2}\vert+\int_0^1\log^{+}\vert(1+e^{2\pi it})/\zeta_{i,2}\vert\mbox{d}t\\
  &=&\log10
 \end{eqnarray*}
 by Jensen's formula \eqref{e: Jensen}. 
 In particular, the entropy is the same for
 $i=1$ and $2$.

 A similar argument for the subgroup $\Lambda=\ll\!\!\e_1,\e_2+\e_3\!\!\gg$
 will show that the entropy is governed by the relation
$$
(1+v_1^{\vphantom2}+v_2^{\vphantom2}\zeta^{-1}_{i,1})
(1+v_1^{\vphantom2}+v_2^{\vphantom2}\zeta^{-1}_{i,2})=0,
$$
 giving entropy
 \begin{eqnarray*}
  \h(\alpha_{M_i}\vert_{\Lambda})&=&
  \log\MM(1+v_1^{\vphantom2}+v_2^{\vphantom2}\zeta^{-1}_{i,1})+
\log\MM(1+v_1^{\vphantom2}+v_2^{\vphantom2}\zeta^{-1}_{i,2})\\
  &=&\log\vert\zeta_{i,1}^{-1}\vert+\int_0^1\log^{+}\vert
\zeta_{i,1}^{\vphantom2}(1+e^{2\pi it})
\vert\mbox{d}t\\
  &&\quad\quad\quad\quad\quad+
  \log\vert\zeta_{i,2}^{-1}\vert+\int_0^1\log^{+}\vert
\zeta_{i,2}^{\vphantom2}(1+e^{2\pi it})
\vert\mbox{d}t\\
  &=&\int_0^1\log^{+}\vert
\zeta_{i,1}(1+e^{2\pi it})
\vert\mbox{d}t+
\int_0^1\log^{+}\vert
\zeta_{i,2}(1+e^{2\pi it})
\vert\mbox{d}t\\
&&\quad\quad\quad\quad\quad\quad\quad\quad\quad\quad\quad\quad\quad\quad
\quad\quad\quad\quad-\log(10),
 \end{eqnarray*}
which is again independent of $i$.

For the general case, let $\Lambda\subset\ZZ^3$ be any rank two
 subgroup of $\ZZ^3$ not already considered.
Find an
element
$\n\in\Lambda\cap(\ZZ^2\times\{0\})$
that is non-zero and
not a non-trivial multiple of any other element of $\Lambda$,
 and
 choose $\m\in\Lambda$ linearly independent to $\n$
with $m_3>0$.
The points $\n$ and $\m$ generate a
 finite-index subgroup of $\Lambda$, so we may assume without
 loss of generality that $\Lambda=\ \ll\!\!\n,\m\!\!\gg$.
Let
 $S=\ZZ[\mathbf u^{\pm\n},\mathbf u^{\pm\m}]$ be the associated
 sub-ring. 
 Consider first the projection $(m_1,m_2,0)$ of $\m$ onto the
$u_1,u_2$-plane: let $w_1=\mathbf u^{\n}$ and
$w_2=u_1^{m_1}u_2^{m_2}$ and write
$S^{\prime}=\ZZ[w_1^{\pm1},w_2^{\pm1}]$. The structure of $M_i$ as
an $S^{\prime}$-module has a single associated principal prime
ideal $\langle\tilde{f}\rangle$ with multiplicity $s({\Lambda})$.
As indicated in the notation, the point of projecting to the
$u_1,u_2$-plane is to remove the variable $u_3$ and ensure that
the sub-module 
structure is the same for $i=1$ and $2$. Now write
$w_3=\mathbf u^{\m}$ so that $S=\ZZ[w_1^{\pm1},w_3^{\pm1}]$, and
notice that $w_3=u_3^{m_3}w_2$. The structure of $M_i$ as an
$S$-module is then determined by the relation
$$
10^{2m_3}\tilde{f}(w^{\vphantom{-m_3}}_1,
w_3^{\vphantom{-m_3}}\zeta_{i,1}^{-m_3})\cdot
\tilde{f}(w_1^{\vphantom{-m_3}},
w_3^{\vphantom{-m_3}}\zeta_{i,2}^{-m_3})=0
$$
in $s(\Lambda)$ copies of the skew plane $\mathbb T^{\Lambda}$.
Since $m_3\neq 0$,
this relation is irreducible over $\mathbb Z$. The entropy of the
action of $\Lambda$ is given by
\begin{multline*}
\frac{1}{s(\Lambda)}\h(\alpha_{M_i}\vert_{\Lambda})\\
=\int_0^1\int_0^1\log\vert
10^{2m_3}\tilde{f}(e^{2\pi it},
e^{2\pi is}\zeta_{i,1}^{-m_3})\cdot
\tilde{f}(e^{2\pi it},
e^{2\pi is}\zeta_{i,2}^{-m_3})
\vert\mbox{d}s\mbox{d}t\\
=2m_3\log10+\int_{0}^{1}\sum\log^+\vert\lambda_s^{\vphantom{m_3}}\zeta_{i,1}^{m_3}\vert
\mbox{d}s+\int_{0}^{1}\sum\log^+\vert\lambda_s^{\vphantom{m_3}}\zeta_{i,2}^{m_3}\vert
\mbox{d}s
\end{multline*}
where the summation is over
the roots of
$\tilde{f}(\lambda_s,e^{2\pi is})=0$.
It follows that the entropy
is independent of $i$.
\end{exam}

\section{Proof of Proposition~\ref{p:meas}}\label{s:proof}

Let $M$ be a Noetherian $R_3$-module, so
$M\cong R_3^k/J$ for some
$R_3$-submodule $J$ of defining relations.
The dual group to $R_3^k$ is
$(\bT^k)^{\ZZ^3}$, so the dual group of
$M$ is the subgroup annihilating $J$.
For $x\in (\bT^k)^{\ZZ^3}$ write
$x_\n^{\vphantom{(1)}}=(x_\n^{(1)},\ldots,x_\n^{(k)})\in\bT^k$ for the coordinate
corresponding to $\n\in\ZZ^3$.

The algebraic $\ZZ^3$-action $\alpha_M$ can be realized as the
usual shift action on the closed, shift-invariant subgroup
of $(\bT^k)^{\ZZ^3}$ defined by
\begin{multline*}
 X_M=\{x\in
 (\bT^k)^{\ZZ^3}\colon f_1(\sigma)(x^{(1)})+\cdots+f_k(\sigma)(x^{(k)})=0\\
 \mbox{ for every }(f_1,\ldots,f_k)\in J\},
\end{multline*}
where $\sigma$ is the $\zd$-shift on $(\bT^k)^{\ZZ^3}$, and
$f(\sigma)$ is the map obtained by substituting the shift into
$f$.

Let $\mathsf X=\mathsf X_M$ be a $\mathbb Z^3$-action with
property $ET$, and write $g(u_3)$ for the expanding polynomial
relation and $f(u_1,u_2)$ for the triangular relation.
The last
paragraph means that $X$ can be thought of as a subshift with
alphabet $\mathbb T^k$ for some finite $k\ge1$.
Define
\begin{equation*}
    S=\{\mathbf m\in\mathbb Z^3\colon m_2\ge0,m_3\ge1\}\cup
\{\mathbf x\in\mathbb Z^3\colon m_2\ge1\}
\end{equation*}
and
\begin{equation*}
    U=\mathbb Z\mathbf e_1.
\end{equation*}
The set $S$ is to be thought of as a lexicographic
`future' for the $\mathbb Z^2$-action generated
by $u_2$ and $u_3$, and the tube $U$ is the `present'.
In our
setting $U$ is a copy of $\mathbb Z$; in the more general setting
of~\cite{Einsiedler02} the set $U$ really is a tube.
We will show later that these notions of `future' and `present'
give rise to the correct entropy.

\begin{figure}[h]
\begin{picture}(200,170)(0,-70)
\put(0,-30){\put(0,0){\vector(4,1){200}}\put(202,48){$u_1$}}
\put(100,-50){\vector(0,1){150}}\put(98,102){$u_3$}
\put(0,-30){\put(0,25){\vector(1,0){220}}\put(222,23){$u_2$}}
\put(-30,-30){\put(130,25){\line(4,1){100}}\put(130,25){\line(-4,-1){100}}}
\multiput(100,-5)(32,8){3}{\put(0,0){\line(0,1){0}}\put(0,30){\line(1,0){0}}\put(30,30){\line(0,-1){30}}
\put(30,0){\line(-1,0){0}}}
\multiput(100,-5)(-32,-8){3}{\put(0,0){\line(0,1){0}}\put(0,30){\line(1,0){30}}\put(30,30){\line(0,-1){30}}
\put(30,0){\line(-1,0){0}}}\put(20,-15){$U$}\put(28,-11){\vector(4,-3){10}}\put(190,70){$S$}
 \thicklines
\multiput(100,25)(32,8){3}{\put(0,0){\line(0,1){60}}\put(0,0){\line(1,0){30}}
\put(30,0){\line(0,-1){80}}}
\multiput(100,25)(-32,-8){3}{\put(0,0){\line(0,1){60}}\put(0,0){\line(1,0){30}}
\put(30,0){\line(0,-1){80}}}
 \put(-150,-90){
\put(100,25){\line(1,0){80}}\put(100,25){\line(4,1){40}}
\put(140,35){\line(4,-1){40}}\put(137,38){$f$}
\put(90,35){\line(0,1){40}}\put(84,53){$g$}}
\end{picture}
\caption{\label{SUpicture}The regions $S$ and $U$ and the shape of
the support of the annihilating relations $f$ and $g$}
\end{figure}
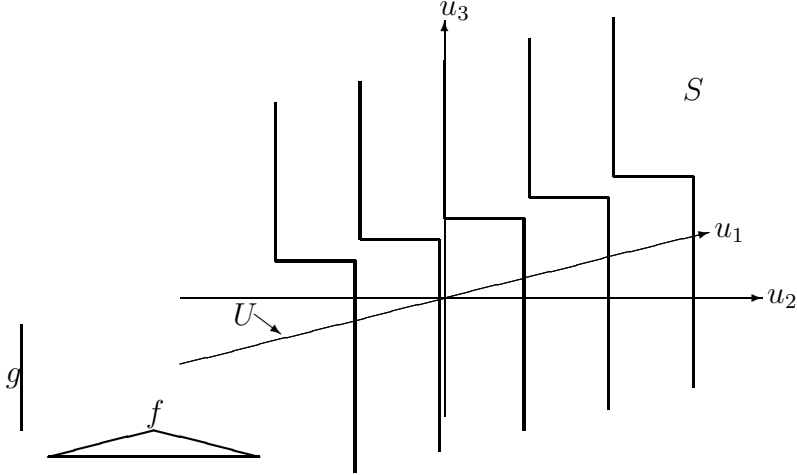

Choose a partition $P$ of the alphabet $\mathbb T^k$ with the
following vertical generating property: for any $\mathbf
n\in\mathbb Z^3$, the $P$-name of the coordinates $x=(x_{\m})$
for $m_1=n_1,m_2=n_2,m_3>n_3$ determines $x_{\mathbf n}$ completely.
Such a
partition exists by the expanding assumption on the relation $g$;
indeed any sufficiently fine partition will do.
This is the only point at which the expanding hypothesis
is essential rather than convenient.
Write $\mathcal P$
for the $\sigma$-algebra defined by the common
refinement $\bigvee_{\mathbf
n\in U}\alpha^{\mathbf n}(P)$.

Let $ [0]=\{x\in X\colon x_{\mathbf n}=0\mbox{ for }\mathbf n\in
S\}. $ There is a natural projection map $\pi$ from $X$ onto
$\left(\mathbb T^k\right)^{S\cup U}$. Let $G=\pi([0])$, and notice
that the properties of $g$ and $f$ ensure that $G$ is a finite
group. Let $\mathcal A$ be the $\sigma$-algebra determined by the
partition $\bigvee_{\mathbf n\in S}\alpha^{\mathbf n}(P)$. Notice
that our assumptions ensure that $\mathcal A$ is the same
$\sigma$-algebra as the pre-image of the whole Borel
$\sigma$-algebra on $\left(\mathbb T^k\right)^{S}$.

Given any point $z\in[0]$, let
\begin{equation*}
    F_{z}(x)=\mu_{x}^{\mathcal A}
    \left([x+z]_{\mathcal A\vee\mathcal
    P}^{\vphantom{\mathcal A}}\right),
\end{equation*}
where $\mu_{x}^{\mathcal A}$ is the conditional measure, and
$[x]_{\mathcal C}$ denotes the $\mathcal C$-atom containing $x$.
Notice that $F_z=F_{z'}$ whenever $\pi(z)=\pi(z')$. If $\zeta\in
G$ and $z\in[0]$ with  $\pi(z)=\zeta$ we set $F_\zeta=F_z$.
Furthermore,
$$
-\log F_0(x)=I_{\mu}(\mathcal P\vert\mathcal A)(x)
$$
is the information function.
More is true: results
from~\cite{Einsiedler02} may be used to show that
\begin{equation}\label{rokhlinmanfredformula}
\int_{X}-\log F_0(x)\mbox{d}\mu=H_{\mu}(\mathcal P\vert\mathcal
A)=h_{\mu} (\alpha\vert_{\Lambda})
\end{equation}
 where $h_{\mu}(\alpha\vert_{\Lambda})$ denotes the entropy with respect
to $\mu$ of the $\mathbb Z^2$-action generated by $\alpha^{\mathbf
e_2}$ and $\alpha^{\mathbf e_3}$. Notice that the second equality
in~\eqref{rokhlinmanfredformula} is a Rokhlin formula, expressing
the entropy of this $\mathbb Z^2$-action as the information
contained in the present given the information contained in the
future with respect to a generator; one subtlety that needs to be
dealt with in proving this is that $\mathcal P$ is not a finite
partition.

We sketch the proof of \eqref{rokhlinmanfredformula}. First, since
the coefficients of the triangular polynomial $f$ are $\pm1$
at the corners, it follows that there exists some non-negative
integer $\ell$ such that
$$\cP\vee\cA=\bigvee_{n=0}^\ell
\alpha^{n\e_1}P\vee\cA.$$
This shows that the integral of the
information function is finite, and that $F_0 (x)>0$ a.e. Furthermore,
again by the properties of the triangular polynomial,
\begin{equation}\label{e: finite}
\bigvee_{\n\in U, |n_1|<N}\alpha^\n P\medspace
\vee\negthinspace\negthickspace\negthickspace\bigvee_ {\n\in
S,|n_1|<N}\alpha^\n P=\bigvee_{\n\in U, 0\leq n_1<\ell}\alpha^\n
P\medspace\vee
\negthinspace\negthickspace\negthickspace\bigvee_{\n\in S,|n_1|<N}\alpha^\n P
\end{equation}
holds for every $N>\ell$. Here we use the fact that the support of $f$ is
triangular. The increasing Martingale theorem for
entropy and~\eqref{e: finite} then
shows~\eqref{rokhlinmanfredformula}.

The next step is to understand how $F_{\xi}$ varies as $\xi$ runs
through the finite group $G$. First, $F_{0}(x)>0$ and
$\sum_{\xi}F_{\xi}(x)=1$ for a.e. $x$. Since $G$ is finite, some
power $m$ of $\alpha^{\mathbf e_1}$ maps $\xi$ to itself, so
$$
F_{\xi}(x)=F_{\xi}(\alpha^{m\mathbf e_1}x).
$$
This implies --- since $\alpha^{\mathbf e_1}$ is totally ergodic
--- that $F_{\xi}(x)=p_{\xi}$ is a.e. constant in $x$. However, for
a different $\xi$ the constant may be different. We claim that
$$
H=\{\zeta\in G\colon p_{\zeta}>0\}
$$
forms a subgroup of $G$. This is proved in~\cite{Einsiedler02} in
detail; to see why it should be true one may argue as follows. For
$\xi_1,\xi_2\in H$, positivity of $F_{\xi_1}(x)$ means that
$x+\xi_1$ is a configuration allowed with positive $\mu$-measure;
positivity of $F_{\xi_2}(x+\xi_1)$ then means that $\xi_1+\xi_2$
is also allowed with the same $\mu$-measure so $\xi_1+\xi_2\in H$.
The full proof requires some care in the removal of a null set of
exceptional behaviour. It follows that
$$
F_{\xi}(x)=\frac{1}{\vert H\vert}\mbox{ for all }\xi\in H.
$$
Notice that the group $H$ is trivial if and only if
$h_{\mu}(\alpha\vert_{\Lambda})$ vanishes.  Since the proposition
is trivial in this case, we assume without loss of generality that
$\vert H\vert=\exp h_{\mu} (\alpha\vert_{\Lambda})>1$. 

At this point we have shown that the measure $\mu$ must be
invariant under translation by a finite, non-trivial, group
(the subgroup $H$ of $G$)
when restricted to a small $\sigma$-algebra of sets (the sets
measurable with respect to $\mathcal A\vee\mathcal P$).
In order
to extend this to invariance on a larger $\sigma$-algebra, we
replicate a version of the argument above on a larger scale.

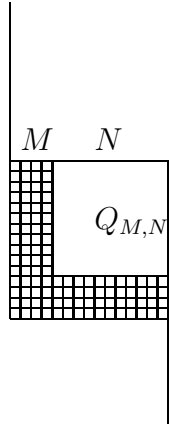
\begin{figure}[h!]
\setlength{\unitlength}{0.07cm}
\begin{picture}(50,100)
\put(0,50){\line(0,1){30}}\put(0,50){\line(1,0){30}}
\put(30,50){\line(0,-1){50}}
\put(0,20){\line(0,1){30}}\put(0,20){\line(1,0){30}}
\put(8,50){\line(0,-1){22}}\put(8,28){\line(1,0){22}}
\multiput(0,20)(0,2){15}{\line(1,0){8}}\multiput(0,20)(0,2){4}{\line(1,0){30}}
\multiput(0,20)(2,0){15}{\line(0,1){8}}\multiput(0,20)(2,0){4}{\line(0,1){30}}
\put(2,52){$M$}\put(16,52){$N$}\put(16,37){$Q_ {M,N}$}
\end{picture}
\caption{\label{MNfigure}The relative size of $B(M,N)$ and $H_N$
as $N\to\infty$, viewed in the plane generated by the last two
standard basis vectors.} 
\end{figure}

Given an integer $N$, define
\begin {eqnarray*}
 S_N&=& \{\mathbf m\in\mathbb Z^3\colon m_2\ge0,m_3\ge N\}\cup
\{\mathbf m\in\mathbb Z^3\colon m_2\ge N\}\\
 U_N&=& \{\mathbf m\in\mathbb Z^3\colon 0\leq m_2 <N, 0\leq
 m_3<N\}
\end {eqnarray*}
and the projection map $\pi_N:X\rightarrow \bigl(\mathbb
T^k\bigr)^{S_N\cup U_N}$ accordingly. If $M<N$, let
$B(M,N)$ denote the set of points $x\in X$ with $x_{\mathbf n}=0$
for all $\mathbf n\in S_N$ except for those coordinates
in the shaded part of
Figure~\ref{MNfigure}. Let $H_N$ be the group defined by the
construction above applied to the scaled action generated by
$\alpha^{N\mathbf e_1},\alpha^{N\mathbf e_2},\alpha^{N\mathbf
e_3}$. Notice that
$$
h_{\mu}(\alpha\vert_{N\Lambda})
=
N^2h_{\mu}(\alpha\vert_{\Lambda}),
$$
so $\vert H_N\vert=N^2\vert H\vert$.
The corresponding
$\sigma$-algebras $\mathcal A_N=\bigvee_{\mathbf n\in
S_N}\alpha^{\mathbf n}(\mathcal P)$ are nested
since $S_N$ moves away from $0$, and
$$
\mathcal
A_N\searrow\mathcal N_X=\{\emptyset, X\}.
$$

Because of the relations $g$ and $f$, there is some fixed number
$K $ such that every little square in Figure~\ref{MNfigure} can at
most contribute a factor
of $K$ to $|\pi_N(B(M,N)|$.
Therefore
$$
\log\vert \pi_N(B(M,N))\cap H_N\vert\le 2MN\log K.
$$
It follows that for any $M$ there is an $N$ such that
$H_N\not\subset B(M,N)$. In particular, if
$Q=Q_{M,n}$ denotes the `inner'
square of side $N-M$ in Figure~\ref{MNfigure} then there exists an
$x\in H_N$ such that $x\vert_Q\neq0$. So we may choose
for each $N$ an element $x^{(N)}\in H_N$ with $x^{(N)}_{\mathbf
n}\neq 0$ for some $\n=\n(N)\in Q_ {M,N} $. After shifting (via
the action) the picture by $\mathbf n$, this gives a point
$y^{(N)}$ with $y_{0}^{(N)}\neq0$. Notice that after shifting the
point $y^{(N)}$ satisfies that $\mu(B+y^{(N)})=\mu(B)$ for
all $B\in\bigvee_{\m\in[-M,M]^3}\alpha^\m P$. By compactness
there is a non-trivial $y\in X$ with the property that $\mu$ is
invariant under translation by $y$ (and therefore by the smallest
closed subgroup of $X$ containing $y$) on the whole
$\sigma$-algebra $\mathcal B_X$.

We are now ready to prove Proposition~\ref{p:meas}.
Let
$$
Y=\{x\in
X\colon\mu\mbox{ is }x\mbox{-invariant}\},
$$
and put $Z=X/Y$. Then 
$(Z,\mu,\alpha)$ has property $ET$ and still satisfies all of the
assumptions of Proposition~\ref{p:meas}. So if
$h(\alpha_Z\vert_{\Lambda})>0$ we would have to find that $\mu$ is
translation-invariant with respect to a non-trivial subgroup,
contradicting the choice of $Y$. It follows that
$h(\alpha_Z\vert_{\Lambda})=0$.

\section{Remarks}

A central question in this type of
algebraic rigidity is the following.
If $\mathsf X$ and $\mathsf Y$ are zero-entropy
mixing algebraic $\mathbb Z^d$-actions, $d>1$, do they exhibit
isomorphism rigidity?
Bhattacharya has shown that the answer is `no'
in general (though he has gone on to show that an
extension of the notion of affine map allows rigidity
to be recovered --- measurable isomorphisms
still arise in rigid families and in particular
are continuous).
Nothing has yet
been shown to preclude a positive
answer to the
following question (cf.~\cite[Conj.~9.1]{klauspims}):
do zero-entropy
mixing algebraic $\mathbb Z^d$-actions, $d>1$ on
connected groups exhibit isomorphism rigidity?

\bibliographystyle{amsplain}
\providecommand{\bysame}{\leavevmode\hbox to3em{\hrulefill}\thinspace}

\end{document}